# A General Model of Ridesharing Services

Carlos F. Daganzo
University of California, Berkeley, CA 94707

Yanfeng Ouyang
University of Illinois, Urbana-Champaign, IL 61820

**Abstract**

The paper presents a general analytic framework to model transit systems that provide door-to-door service. The model includes as special cases *non-shared taxi* and *demand responsive transportation* (*DRT*). In the latter we include both, paratransit services such as *dial-a-ride* (*DAR*), and the form of ridesharing (*shared taxi*) currently being used by crowd-sourced taxi companies like Lyft and Uber. The framework yields somewhat optimistic results because, among other things, it is deterministic and does not track vehicles across space. By virtue of its simplicity however, the framework yields approximate closed form formulas for many cases of interest.

## 1. The model

The model pertains to a region of size $R$ [km$^2$] with a steady and uniform many-to-many demand density i.e, where the trip origins and destinations are uniformly distributed in space and time. The trip origination rate is denoted $\lambda$ [pax/hr-km$^2$]. This demand is served by a fleet of $m$ vehicles. Each of these vehicles can travel at speed $v$ [km/hr]. Passengers are assumed to board and alight quickly, so the times for these actions are neglected. A central controller uses an algorithm to assign callers to vehicles and route the vehicles. The algorithm can create a buffer of unassigned callers to make assignments more efficient. The size of this buffer, $n_w$, is a control variable. The algorithm can also artificially restrict the maximum number of occupants in a vehicle to any number $c$ that does not exceed the vehicle's capacity. This number is a second control variable. The model about to be described has the ability to represent different algorithms. These can range from individual taxi service to dial-a-ride, and include the type of ride-sharing used by crowd-source taxi companies.

In the model, the current status of the transit system is characterized by a state vector $\boldsymbol{n}$ that is tracked over time with a set of differential dynamic equations. This vector is composed of the number of vehicles under different workloads, but without reference to their spatial position. A vehicle's workload is only characterized by a tuple of non-negative integers $(i, j)$: the first index is the number of passengers in the vehicle and the second the number of passengers assigned to it for pick up. We shall use $n_{ij}$ to denote the number of vehicles under workload level $(i, j)$, so $\boldsymbol{n} = \{n_{ij}\}$. This type of space-less, deterministic model is only approximate but can be useful for policy analysis.

A vehicle can change workload level in three ways:

(i) An assignment: $(i, j)$ changes to $(i, j+1)$. The rate at which this occurs is denoted $a_{ij}$ [veh/hr].
(ii) A pickup: $(i, j)$ changes to $(i+1, j-1)$. The rate at which this occurs is denoted $p_{ij}$ [veh/hr].
(iii) A delivery: $(i, j)$ changes to $(i-1, j)$. The rate at which this occurs is denoted $d_{ij}$ [veh/hr].

Figure 1 encapsulates all this information in the form of a network. The nodes are workload levels and the links transitions. Double arrows pointing down are assignments, slanted arrows pointing up collections and dashed arrows pointing left deliveries. Note that the network has $N = (c+1)(c+2)/2$ nodes and $A = P = D = c(c+1)/2$ links of each type.

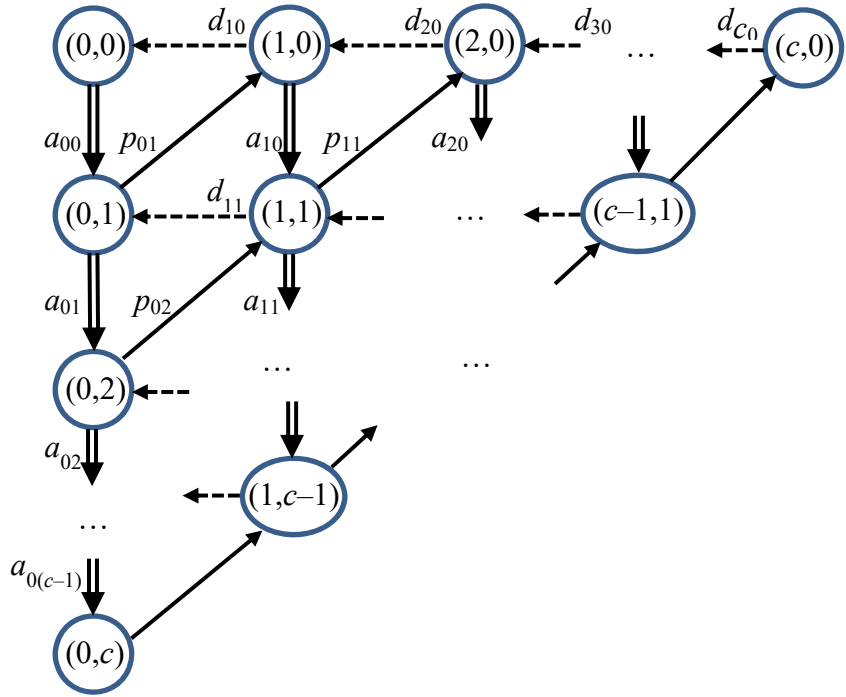

Figure 1. Workload transition network.

The control algorithm determines how the vectors, $a=\{a_{ij}\}$, $p=\{p_{ij}\}$, $d=\{d_{ij}\}$ depend on $n$; i.e., the vector functions $\{a(n), p(n), d(n)\}$. If we use row vector notation and $A$, $P$ and $D$ for the link-node incidence matrices of each link type we can write the system's dynamic equations as:

$$d\mathbf{n}/dt = a(n)A + p(n)P + d(n)D. \quad (1)$$

Clearly, a steady state solution $n$ must satisfy the above with $d\mathbf{n}/dt = \mathbf{0}$; i.e., the following flow conservation equations:

$$a(n)A + p(n)P + d(n)D = \mathbf{0}. \quad (2)$$

Note that one of the scalar equations in (2) is always redundant -- as must occur for network flow problems. Thus, (2) involves $N-1$ independent equations. Since $n$ has dimension $N$, there is one more variable than there are equations. However, if $m$ is given $n$ must also satisfy an additional equation expressing the fleet size constraint. If we use a boldface $\mathbf{1}$ for a column vector of 1's, this constraint is:

$$\mathbf{n1} = m. \quad (3)$$

Equations (2) and (3) together involve $N$ equations for $N$ unknowns, suggesting that the problem may have a unique solution. The equations are non-linear, however, and of course, $n$ must also be non-negative. Thus, a steady state solution may not exist. This may happen for example if $m$ is too small to serve the demand.

The equations can be used to assess the performance of any algorithm for any given level of demand. This evaluation can be done generically for all system configurations $(\lambda, R, v, m)$ with only two degrees of freedom because we can choose the units for distance and time so that $R = 1$ and $v = 1$. The degrees of freedom are then the fleet size $m$ and the rescaled demand density in the new intrinsic units, which we denote $\pi \equiv \lambda R^{3/2}/v$. Think of the latter as the number of calls that arrive in the time it takes a vehicle to cross the region.

Treating the problem in this intrinsic system of units is useful. It reveals for example that the minimum fleet size that is required to serve the demand with a given algorithm be a univariate function: $m_c = m_c(\pi)$. The solution can always be expressed in an arbitrary system of units by rescaling the solution *a posteriori*; e.g., by writing $m_c = m_c(\lambda R^{3/2}/v)$. Similarly, any other measure must be of the form $f(\pi, m)$ in intrinsic units, and these solutions can also be rescaled to arbitrary units. Therefore, the rest of this paper will use intrinsic units. Sections (1)-(3) derive the functions $m_c(\pi)$ and $f(\pi, m)$ for three algorithms that emphasize three different objectives: (i) level of service; (ii) cost reduction, and (iii) passenger capture.

## 2. Level of service emphasis (non-shared taxi)

Level of service is emphasized by setting $n_w = 0$ and $c = 1$. This ensures that people can be told upon calling when they will be collected and delivered. Furthermore, to reduce the callers' waits the algorithm assigns to every call the closest available taxi. The preceding stipulations fully define the operation, since by setting $c = 1$ the vehicle operator has no choices to make regarding delivery. In essence the operation becomes that of a taxi service without ridesharing. This is the taxi operating scheme described in Section 7.1 of Daganzo and Ouyang (2018).

To analyze the scheme, consider Figure 2, which depicts the algorithm's network diagram. Because the network has only three links, the flow conservation equation (2) reduces to:

$$a_{00}(\mathbf{n}) = p_{01}(\mathbf{n}) = d_{10}(\mathbf{n}). \tag{4}$$

The three members of (4) can be expressed approximately using the formula for the expected distance between a random point in a region (with $R=1$) and the closest of $r$ random points, $\delta(r)$. This approximation is optimistic as it assumes that assignments are made instantaneously, whereas in reality there is a delay.[1] For consistency, the approximation is also made when analyzing the other service types. The formula is (see for example, Daganzo and Ouyang, 2018, for a proof):

$$\delta(r) \approx k r^{-1/2}, \tag{5}$$

where $k$ is a dimensionless constant that depends on the network topology. It is $k \approx 0.63$ for networks that resemble rectangular grids. Equation (5) is asymptotically exact when $n$ approaches infinity. It is also a very good approximation for low values of $n$ if the region in question is fairly round in shape, e.g. resembling a square or circle.

---

[1] Delay in assignment leads to the relatively rare (but real) scenario in which the same available vehicle is the closest option for multiple near-simultaneous calls.

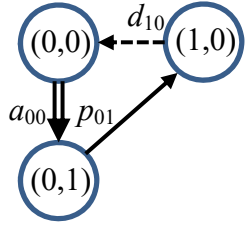

Figure 2. Workload transition network under level of service emphasis.

Equation (5) means that for regions of these shapes the average distances (and times since $v=1$) traveled by a taxi to deliver and pick up a customer can be respectively approximated by $\delta(1) \approx k$ and $\delta(n_{00}) \approx k(n_{00})^{-1/2}$. Keeping this in mind we see from Little's formula (applied to each of the nodes) that our steady state flows can be approximated by the following:

$$a_{00}(\mathbf{n}) = \pi, \tag{6a}$$
$$p_{01}(\mathbf{n}) \approx n_{01}/\delta(n_{00}) \approx n_{01}(n_{00})^{1/2}/k, \tag{6b}$$
$$d_{10}(\mathbf{n}) \approx n_{10}/\delta(1) \approx n_{10}/k. \tag{6c}$$

Thus, (4) becomes:

$$\pi \approx n_{01}(n_{00})^{1/2}/k \approx n_{10}/k. \tag{7}$$

Note that these equations are non-linear only because of the factor $(n_{00})^{1/2}$ that arises from the function $\delta$ in (6b). Therefore, we now simplify them by conditioning on $n_{00}$. Since the flow conservation equations will also be simplified in the next sections with different conditioning variables, it will be convenient to denote the conditioning variable by $n$. In the present case $n \equiv n_{00}$ and the flow conservation equations (7) become:[2]

$$n_{00} = n, \tag{8a}$$
$$n_{01} \approx k\pi/n^{1/2}, \tag{8b}$$
$$n_{10} \approx k\pi. \tag{8c}$$

The fleet size is then

$$m \approx n + k\pi n^{-1/2} + k\pi. \tag{9}$$

---

[2] Note that $n$ is the size of the choice set considered by the algorithm when making an assignment for pickup. This will also be the case in the following sections.

This is a convex function, which reaches a unique minimum—the critical fleet size:

$$m_c \approx 3(k\pi/2)^{2/3} + k\pi. \qquad (10)$$

We can also find expressions for any other measure of performance $f$ that is a function of $\mathbf{n}$; e.g., for the fraction of time $f_i$ that taxis spend idling, which is $f_i = n_{00}/m$ as per Little's formula. Another example is the ratio $f_t$ of the average passenger door-to-door travel time (including waiting), which is $(1-n_{00}/m)(m/\pi)$, to the time it takes to drive directly from the origin to the destination $\delta(1) \approx k$; i.e:

$$f_t = (m - n_{00})/(k\pi). \qquad (11)$$

In the following we shall focus on $\{f_t;\, m\}$ because this curve captures very well the key tradeoff between user time and the cost of providing service—since the latter (externalities included) is roughly proportional to $m$.

This measure of performance (and any others) can be plotted as curves in the $(f, m)$ plane. These curves can also be expressed in parametric form with $n$ as the parameter by inserting (8) into the formulae for $f$. In the case of (11) the expression is $\{(m-n)/(k\pi);\, n + k\pi/n^{1/2} + k\pi : n > 0\}$. The instance with $\pi = 100$ and $k = 0.63$ is plotted in Figure 3. Note that the pictured curve has two branches. The upper branch is obviously inefficient and should not arise. For this reason it is shown by means of a dotted line. The solid part shows how door-to-door travel time declines with fleet size (cost) in the range $m > m_c$.

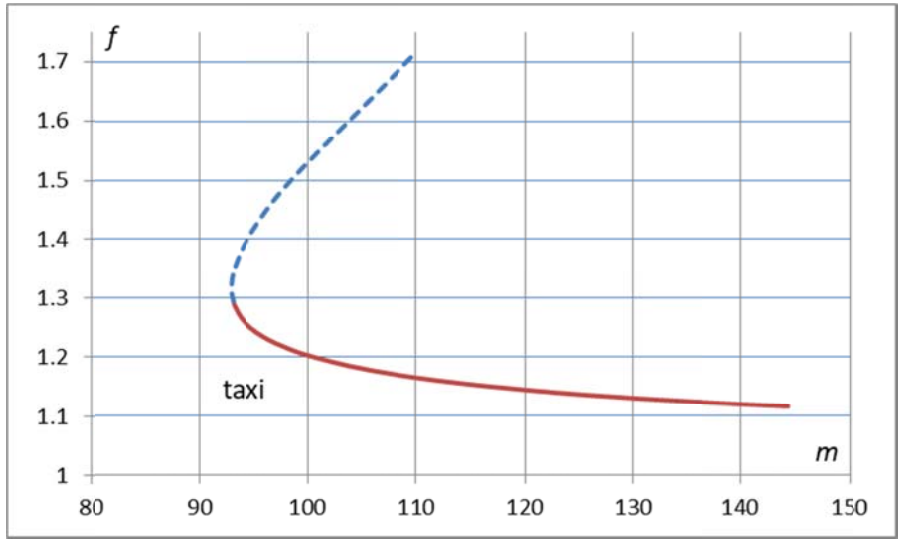

Figure 3. User travel time vs. fleet size (proxy for system cost).

## 3. Cost reduction emphasis (dial-a-ride)

Operating cost is emphasized by maximizing vehicle productivity. This is achieved by setting the number of unassigned callers at "home," $n_w$, equal to some large value, and setting $c$ equal to the vehicles' capacity. To further increase productivity, passenger occupancies of all circulating vehicles are kept at the maximum allowed level. Thus, upon delivering a passenger the driver then picks up a passenger, and chooses the nearest of the $n_w$ unassigned passengers. And after a pickup, the driver always chooses for delivery the vehicle's occupant with the closest destination, regardless of how long other occupants have been onboard. This efficient (but drastic) operating scheme is intended to mimic paratransit demand responsive (dial-a-ride) services. It was analyzed in Daganzo (1977) and further elaborated in Section 7.5 of Daganzo and Ouyang (2018).

Figure 4 shows the network diagram for this form of operation. It only involves three nodes. As we shall see, it can also be solved in closed form.

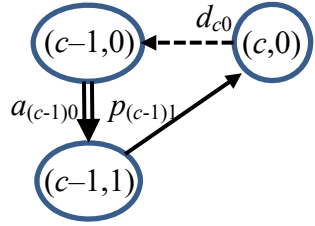

Figure 4. Workload transition network under cost reduction emphasis.

First we write expressions for the link flows. Note that the size of the choice set for pickup is $n_w$ and the size for delivery is $c$. Therefore we choose $n \equiv n_w$ and the counterparts of (6) become:

$$a_{(c-1)0}(\bm{n}) = \pi, \tag{12a}$$
$$p_{(c-1)1}(\bm{n}) \approx n_{(c-1)1}/\delta(n) \approx n_{(c-1)1}(n)^{1/2}/k, \tag{12b}$$
$$d_{c0}(\bm{n}) \approx n_{c0}/\delta(c) \approx n_{c0}(c)^{1/2}/k. \tag{12c}$$

Next, we use the flow conservation equations to express the vector $\bm{n}$ in terms of $n$. Recall that all vehicles are immediately assigned a passenger when the occupancy is less than $c$. As a result vehicles spend no time with workload $(c-1, 0)$ and therefore $n_{(c-1)0} = 0$. Moreover, the counterpart of (7) is now $\pi \approx n_{(c-1)1}(n)^{1/2}/k \approx n_{c0}(c)^{1/2}/k$. Thus, we now have:

$$n_{(c-1)0} = 0, \tag{13a}$$
$$n_{(c-1)1} \approx k\pi n^{-1/2}, \tag{13b}$$
$$n_{c0} \approx k\pi c^{-1/2}. \tag{13c}$$

Therefore, the fleet size is:

$$m \approx k\pi n^{-1/2} + k\pi c^{-1/2}. \tag{14}$$

Logically, this is minimized by setting $n = \infty$, so:

$$m_c \approx k\pi c^{-1/2}. \tag{15}$$

This is, of course, utopian because passenger waits at home would be infinite.

Note that the average passenger's door-to-door travel time is the ratio of the number of passengers in the system to the passenger generation rate $\pi$. The numerator of this ratio is $(n+mc)$ since it must combine the $n$ unassigned passengers waiting at home and the $c$ passengers who are either riding or assigned to each of the $m$ vehicles. If people were to drive, their travel time would be $k$. Thus, the ratio of the average door-to-door travel time to the driving time is $f_t = (n+mc)/(k\pi)$. In terms of $n$ it is:

$$f_t = n/(k\pi) + cn^{-1/2} + c^{1/2}. \tag{16}$$

This relation and (14) allow us to write the curve $\{f_t; m\}$ in parametric form. It is depicted in Figure 5 for the case with $\pi = 100$ and $k = 0.63$ for $c = 2, 3$ and $5$ using the same conventions as in the previous section. Similar to the taxi case, for each $c$, the plotted curve has two branches. The branch with larger $m$ is obviously inefficient and would never be chosen – it corresponds to small values of the control variable $n \equiv n_w$, which result in unnecessarily long in-vehicle riding times. Furthermore, the model is a good approximation to real-world operations only when the number of waiting passengers, $n$, is considerably greater than 1, such that randomness in the passenger arrivals does not significantly impact system operations. This imposes an upper bound on the fleet size in light of (14). For example, if we want to ensure $n \geq 2$, then the upper bound of fleet size is

$$\hat{m} = 2^{-1/2}k\pi + k\pi c^{-1/2}.$$

The figure also includes the taxi curve for comparison purposes. Note how the dial-a-ride algorithm can function with a smaller fleet size, and that there is a gap between both sets of curves. It stands to reason that it should be possible to fill this gap with other algorithms. The next section explores this idea.

4. **Passenger capture emphasis (shared taxi)**

Our third scenario is intended to resemble the ridesharing systems provided by for-profit taxi companies. We will thus, assume that $c > 1$ as in Section 3. However, because these companies operate in a highly competitive environment they cannot afford to hold reservoirs of unassigned callers to improve efficiency (otherwise waiting customers could switch to competitors) we shall assume that $n_w = 0$, as in Section 2. Thus the scenario about to be studied is a hybrid of taxi and dial-a-ride, which should provide an in-between level of service.

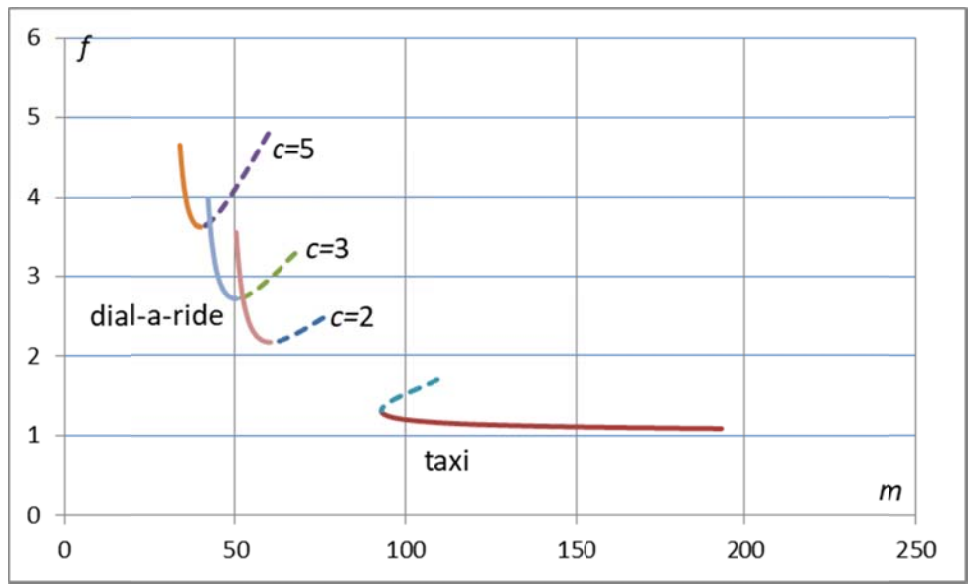

Figure 5. User travel time vs. fleet size (proxy for system cost).

To maximize passenger capture we shall assume that pickups take priority over deliveries and that callers are assigned to the closest of the $n$ available vehicles. The set of available vehicles can be defined in different ways, however. For example, it can include all vehicles with $(i+j) < c$ so that

$$n = \sum_{i+j<c} n_{ij}. \tag{17a}$$

Or it can include only those taxis with $i = 0$ and $j < c$ so that

$$n = \sum_{j<c} n_{0j}. \tag{17b}$$

The first option results in quicker pickups, while the second in more predictable rides—since it prevents taxis from initiating pickups when there are passengers on board. Proprietary algorithms probably do something in between these two extremes, so by studying both once can estimate what real companies do. Figures 6(a) and 6(b) include the network diagrams for these two types of operation with $c = 2$. Note that the dashed delivery arrows only exist when vehicles have no outstanding pickup assignments. In the following, the procedure is applied to the system corresponding to Figure 6(b).

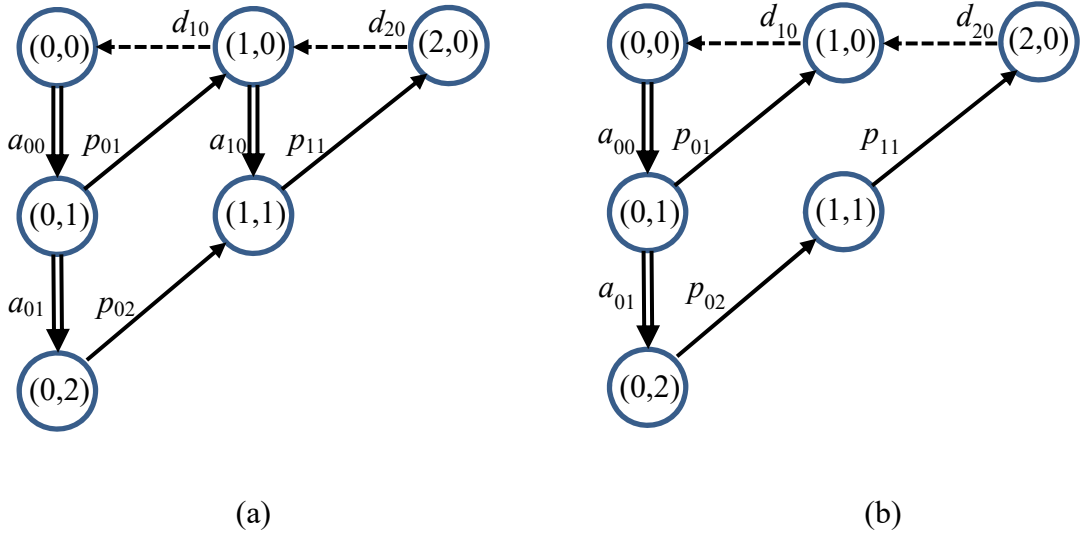

(a)　　　　　　　　　　　　　　　　(b)

Figure 6. Workload transition network under passenger capture emphasis and two call assignment protocols: (a) to the nearest vehicle with room; (b) to the nearest vehicle with room and no onboard passengers.

In this case (17b) reduces to:

$$n = n_{00} + n_{01}; \tag{18}$$

and, given our routing assumptions, the link flow functions are:

$$a_{ij}(\mathbf{n}) = \pi\, n_{ij}/n, \qquad (ij) = (00), (01) \qquad (19a)$$
$$p_{ij}(\mathbf{n}) \approx n_{ij}/\delta(n) \approx n_{ij}(n)^{1/2}/k, \qquad (ij) = (01), (11), (02) \qquad (19b)$$
$$d_{ij}(\mathbf{n}) \approx n_{ij}/\delta(i) \approx n_{ij}(i)^{1/2}/k. \qquad (ij) = (10), (20) \qquad (19c)$$

Flow conservation at the six network nodes yields the following five independent equalities:

$$\pi n_{01}/n = n_{02}(n)^{1/2}/k = n_{11}(n)^{1/2}/k = n_{20}(2)^{1/2}/k \quad ; \text{ for nodes } (02), (11), (20) \qquad (20a)$$
$$\pi n_{00}/n = n_{10}/k \qquad\qquad ; \text{ for node } (00) \qquad (20a)$$
$$n_{01}(n)^{1/2}/k + n_{20}(2)^{1/2}/k = n_{10}/k \qquad ; \text{ for node } (10). \qquad (20a)$$

Now treat $n$ as a constant and solve the six linear equations (18) and (20) for the six state variables $n_{ij}$ in terms of $n$. Because this is straight-forward the steps are omitted. The final results are:

$$n_{00} = n\left[k\pi + n^{\frac{3}{2}}\right]\Big/\left[2k\pi + n^{\frac{3}{2}}\right],\ n_{01} = k\pi n\Big/\left[2k\pi + n^{\frac{3}{2}}\right],\ n_{02} = n_{11} = (k\pi)^2\Big/\left[2k\pi n^{\frac{1}{2}} + n^2\right],$$
$$n_{10} = k\pi\left[k\pi + n^{\frac{3}{2}}\right]\Big/\left[2k\pi + n^{\frac{3}{2}}\right],\ n_{20} = (k\pi)^2\Big/\left[2\sqrt{2}k\pi + \sqrt{2}n^{\frac{3}{2}}\right]. \qquad (21)$$

We can now express in terms of $n$ any measure of performance that combines the state variables. For example, the user time (as compared to driving) vs. societal cost that was analyzed in previous sections is now:

$$\{f_t; m\} = \left\{\frac{1}{k\pi}\sum_{ij}(i+j)n_{ij}(n); \sum_{ij} n_{ij}(n) : n > 0\right\}. \qquad (22)$$

The minimum fleet size is the value of $m$ where curve (22) bends back on itself. It is given analytically by:

$$m_c = \inf_{n>0}\left\{\sum_{ij} n_{ij}(n)\right\}. \qquad (23)$$

As in previous sections curve (22) is plotted for $\pi = 100$ and $k = 0.63$ on the same plane with the other curves; see Figure 7. This figure also includes a curve, labeled "(a)" with the result for the algorithm in part (a) of Figure 6. The formulas underpinning this latter result were obtained with the same method we have been using. They are not given here in the interest of brevity. Note how the critical fleet sizes for the two ride-sharing algorithms lie in between those of non-shared taxi and dial-a-ride, and how their performance curves fill the gap between the two older technologies—as we had anticipated.

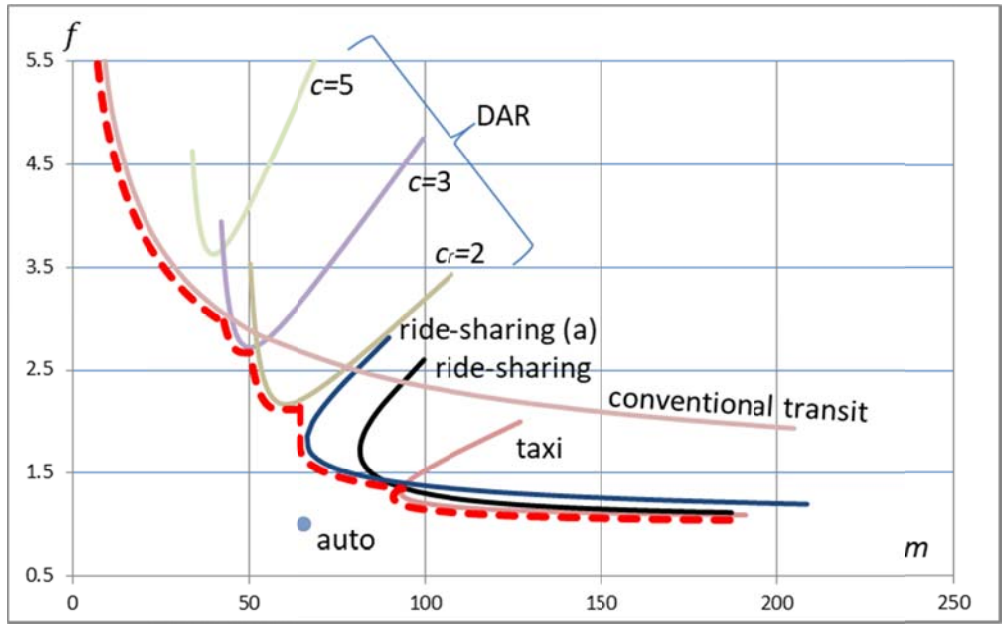

Figure 7. User travel time vs. fleet size (proxy for system cost) for $\pi = 100$.

Figure 7 also includes a solid dot representing travel by private automobile. The abscissa is the number of cars in circulation, assuming a square region with a dense street grid. Since the average trip length is $(2/3)R^{1/2}$, Little's formula yields this number as: $(2/3)\lambda R^{3/2}v^{-1} \equiv (2/3)\pi \approx 66$.

The figure also displays a curve representing the performance of conventional transit. This curve is based on the following assumptions: (i) routes form a uniform square grid over a square region; (ii) the $m$ vehicles are evenly distributed across all routes; (iii) transfer and station dwell times are negligible; and (iv) the passengers' walking speed is 1/10 of the vehicle speed $v$. The logic is as follows. When the route spacing is $S$, there are a total of $2R^{1/2}/S$ bus routes, each requiring buses to spend $2R^{1/2} v^{-1}$ time units to perform a round trip. Hence, the number of vehicles serving each route is $mS/(2R^{1/2})$, and the headway is $4R/(vmS)$. Now, if passengers wait for half of the headway at the boarding station, and walk a total of $S$ distance units to access stations at both ends of the trip, then the expected waiting time is $2R/(vmS)$. Furthermore, the expected walking time is $10S/v$. Finally, noting the average riding time is still $(2/3) R^{1/2} v^{-1}$ we see that the expected door-to-door time including all components is

$$10S/v + 2R/(vmS) + (2/3) R^{1/2} v^{-1} \geq (2/3) R^{1/2} v^{-1} + 4 [5R/m]^{1/2} v^{-1}.$$

The inequality holds because its left side is an EOQ (economic order quantity) formula with respect to $S$, and the right side is the EOQ's minimum. The right side is the door-to-door time when $S$ is chosen optimally. By setting $R=v=1$ we obtain our dimensionless estimate of door-to-door travel time: $f = 1+6[5/m]^{1/2}$. This is the formula of the plotted curve.

Note that, unlike the corresponding curves for the flexible transit modes, the curve for conventional transit decreases monotonically with increasing $m$ and approaches the performance of the automobile ($f = 1$) for very large $m$. This occurs because with conventional transit people are asked to walk to designated stations, where they can be efficiently served in large groups.

The thick dashed curve bordering the bottoms of the transit curves is the "Pareto frontier" of all the transit modes; i.e., the set of points whose cost-time coordinates are not improved by any other point on any curve. The section of the frontier that overlaps with a particular modal curve defines the niche of that mode.

For comparison purposes, Figure 8 displays the same diagrams for two scenarios with higher demand: $\pi = 1{,}000$ and $\pi = 10{,}000$. Note by comparing the different figures how the fleet size needed to achieve a certain door-to-door travel time standard increases with demand. Note as well how the niche of conventional transit increases with demand, but that of flexible transit decreases. Also note how the auto dots are closer to the Pareto frontier. All this suggests that the bigger and denser a city the more favorable is its environment for conventional transit.

## 5. Stochastic effects: a comparison with simulated results

The results from Sections 2-4 are deterministic in the sense that they ignore the stochastic fluctuations of the state variables that arise in the real world. This simplification should cause the deterministic formulae to under-predict both travel times and minimum required fleet sizes. These under-predictions should be particularly large when the stochastic fluctuations are large compared with the means; i.e., when these means are low. This section quantifies the under predictions by juxtaposing simulated results with the analytical curves. The resulting diagrams show that the analytical formulas provide qualitatively correct comparisons across the different service types.

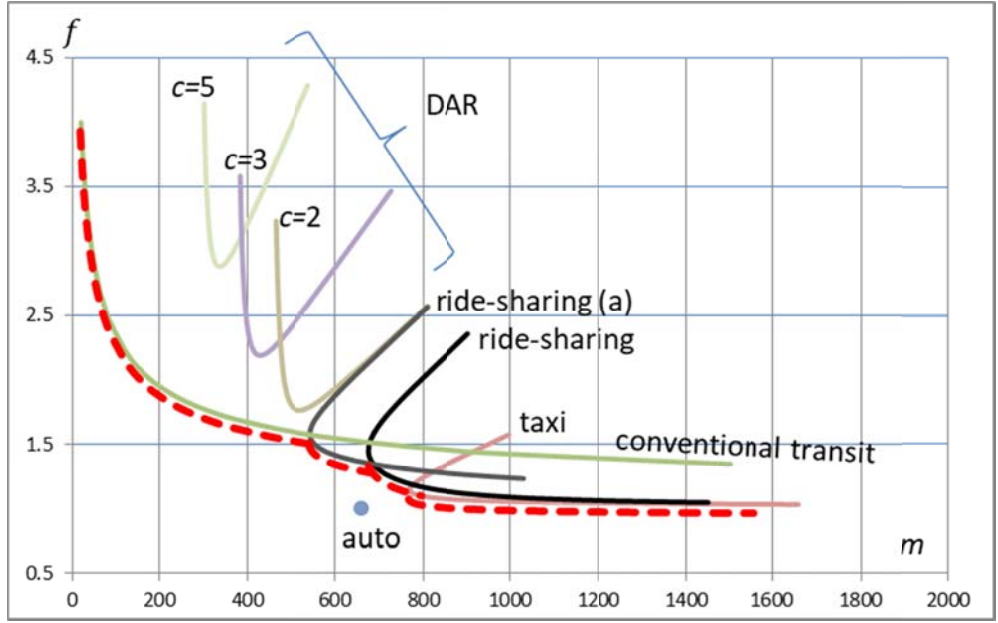

(a) $\pi=1000$

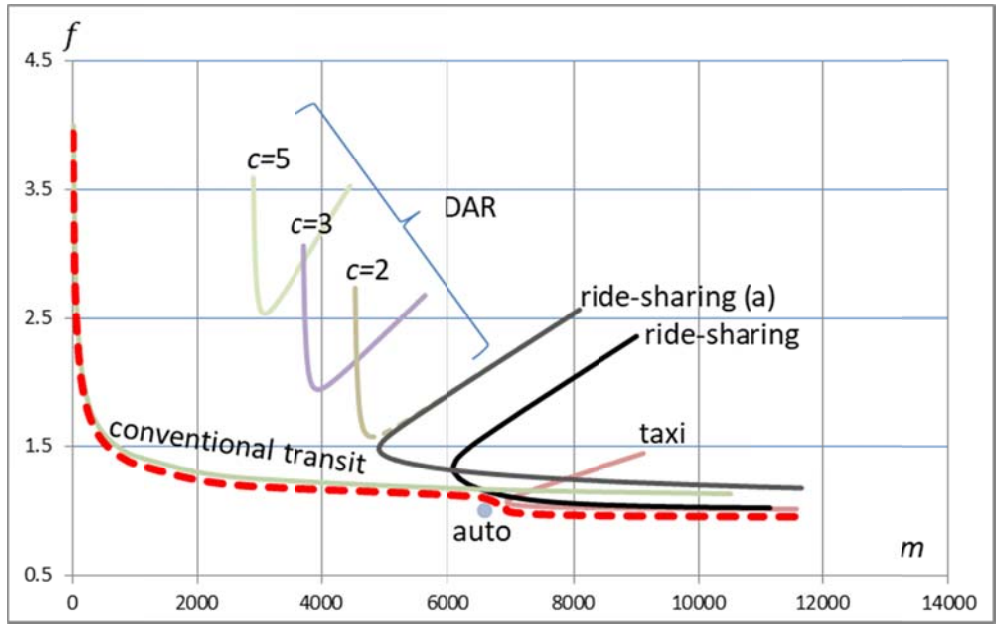

(b) $\pi=10,000$

Figure 8. User travel time vs. fleet size (proxy for system cost) under higher demand.

An agent-based simulation program was written that tracked both, passengers and vehicles in space and time as they followed the rules of our algorithms under steady-state Poisson demand. After a "warm up" period including 500 passengers, each simulation run recorded the travel times of the following 10,000 passengers, and their average was compared against that predicted by the analytical formulae. Multiple runs with different fleet sizes were completed for each of the scenarios of Figures 7 and 8.

Since stochastic fluctuations play the largest role when the demand is low, we only present below the results for the worst case ($\pi$=100) which was depicted in Figure 7. They are displayed in Figure 9, which also includes the corresponding deterministic curves. Triangles, diamonds and squares are used to differentiate the different algorithms; i.e., the results for taxi, shared-taxi (type a) and shared taxi (type b) in part (a) of the figure and the three cases of DAR, with $c = 2, 3$ and 5, in part (b).

In Figure 9(a) results were produced only for fleet sizes able to sustain a steady state. This is why the collections of triangles dots and squares do not continue farther to the left. Also note how in this figure the analytical curves lie below the corresponding symbols—this is the abovementioned under-prediction. Nonetheless, the symbols display very similar relative positions as the curves, and tell the same qualitative story. Note for example how the lower envelope of the three curves defined by the symbols starts with shared-taxi-"a" on the left, then there is a small range with shared-taxi-"b", and the envelope ends with ordinary taxi on the right. This is the same pattern described by the analytic curves.

Note as well how the vertical displacements in Figure 9(a) are larger for small $m$. This should not be surprising, since in the case of taxi the errors are mostly driven by fluctuations in the number of idle taxis. Since for a given demand $\pi$ the average number of idle taxis decreases as the fleet size $m$ decreases, stochastic fluctuations should then matter more. For example, when the expected number of idle taxis is comparable with one, fluctuations in the number idle taxis will frequently drive the number of idle taxis to zero, forcing some passengers to wait for an assignment. This passenger delay increases travel time in a way that is not captured by the analytical model. Fortunately, because the curves are steeper when their under-predictions are larger, we find that the analytical predictions can be improved markedly if the curves are shifted to the right by about 20 vehicles; i.e., one pretends that the actual fleet size is reduced by 20 vehicles. With this correction, the analytical model is fairly accurate and can be used to make rough but reasonable quantitative predictions. Note in particular how the minimum fleet sizes predicted by the simulation {90, 100, 110} and the model {67, 82, 93} differ by about 20 vehicles.

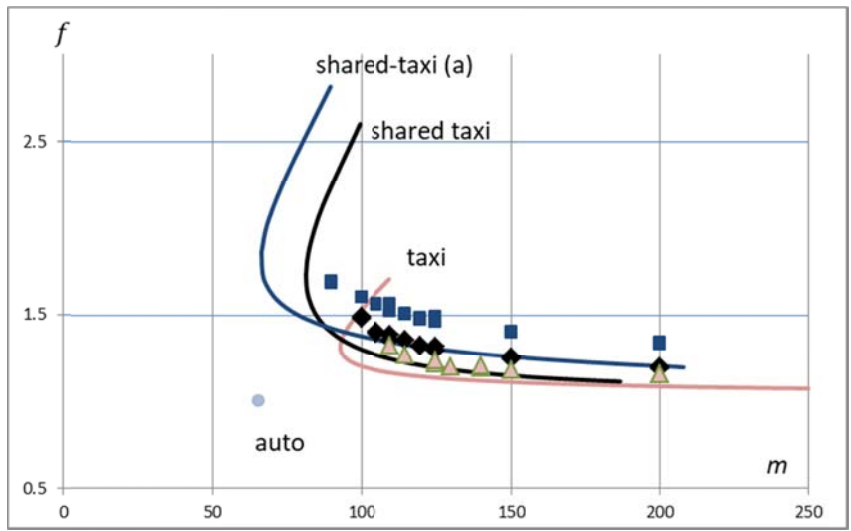

(a) Shared taxi and taxi

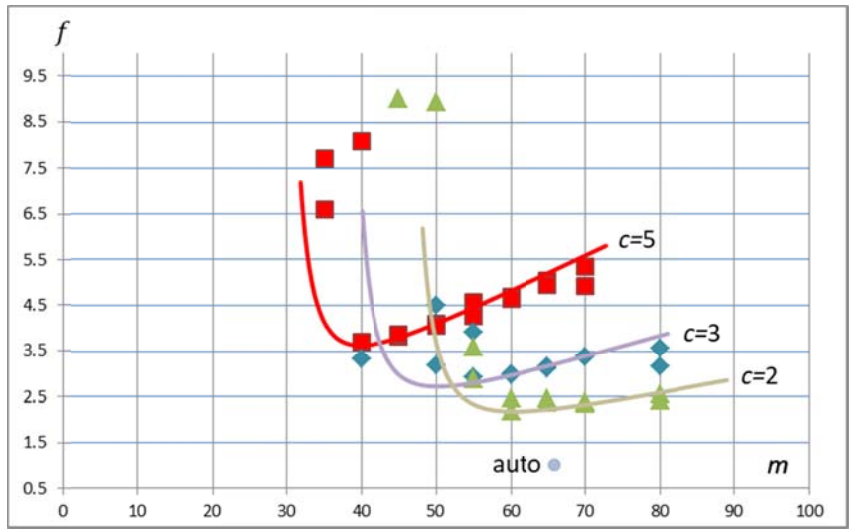

(b) DAR

Figure 9. Simulation of various ride-sharing strategies.

For DAR the results are somewhat better; see Figure 9(b). This occurs because in this case the only state variable that fluctuates randomly is the number of passengers waiting at home, which is usually a large number. Thus, the stochastic fluctuations do not bias the results as much. As a result the symbols are closer to the curves and the simulated minimum fleet sizes {45, 60, 60}, which include considerable statistical error, are also closer to the model predictions {40, 50, 61}.

## 6. Concluding Remarks

The results in this note are just a beginning. Although the provided method is reasonably accurate, especially for large systems with heavy demand, it is not perfectly realistic. The main simplification is that the assignment algorithms used by service providers are likely to be much more complicated than what we studied. For example, the analysis assumes that passenger assignments in the taxi cases are instantaneous, but in reality there is a delay. Furthermore, service providers may consider the alignment of destinations in making assignment decisions, and not just origins. Despite all these shortcomings, however, the evidence seems to suggest that the analytical model can be used by taxi companies and government agencies to systematically explore operating and pricing strategies. They could also be used by cities to begin to understand how taxi companies might respond to different forms of regulation.